\documentclass[preliminary,copyright,creativecommons]{eptcs}
\pdfoutput=1
\providecommand{\event}{ACT 2023} %
\usepackage{breakurl}             %
\usepackage{underscore}           %
\usepackage{macros}

\usepackage[leftmargin=1em,rightmargin=1ex,vskip=1ex]{quoting}
\usepackage{enumitem}
\usepackage{lscape}
\usepackage{unicode-mario}
\usepackage{quiver}
\usepackage{framed}

\title{Collages of String Diagrams}
\author{Dylan Braithwaite
  \institute{University of Strathclyde}
  \email{dylan.braithwaite@strath.ac.uk }
  \and
  Mario Rom\'an
  \institute{Tallinn University of Technology}
  \email{mroman@ttu.ee}
}

\def\titlerunning{Collages of String Diagrams}
\def\authorrunning{Braithwaite and Rom\'an}

\makeatletter
\def\@copyrightspace{\relax}
\makeatother

\begin{document}
\maketitle

\begin{abstract}
  We introduce collages of string diagrams as a diagrammatic syntax for gluing multiple monoidal categories.
  Collages of string diagrams are interpreted as pointed bimodular profunctors.
  As the main examples of this technique, we introduce string diagrams for bimodular categories, string diagrams for functor boxes, and string diagrams for internal diagrams.
\end{abstract}

\section{Introduction}

String diagrams are a convenient and intuitive, sound and complete syntax for monoidal categories \cite{joyal91}. 
Monoidal categories are algebras of processes composing in parallel and sequentially \cite{macLane1971}; string diagrams formalize the process diagrams of engineering \cite{boisseau21,bonchi19}. Formalization is not only of conceptual interest: it means we can sharpen  our reasoning, scale our diagrams, or explain them to a computer \cite{patterson2021}.  

However, the formal syntax of monoidal categories is not enough for all applications and, sometimes, we need to extend it. Functor boxes allow us to reason about translations between theories of processes \cite{cockett1999linearly,mellies06}, ownership \cite{nester20ledger}, higher-order processes \cite{alvarezpicallo2021}, or programming effects \cite{pirog16:stringsfreemonads}. Quantum combs not only model some classes of supermaps \cite{chiribella2009,coeckeFS16,hefford2022coend}, but they coincide with the monoidal lenses of functional programming \cite{boisseauG18,clarke2022profunctor,spivak2022generalized} and compositional game theory \cite{compositionalGameTheory,bolt:bayesian}. Premonoidal categories, which appear in Moggi's semantics of programming effects \cite{moggi91,levy2004,staton13}, are now within the realm of string diagrammatic reasoning \cite{roman:effectful}. Internal diagrams extend the syntax of monoidal categories allowing us to draw diagrams inside tubular cobordisms and reason about topological quantum field theories \cite{bartlett2015modular}, but also coends \cite{roman21} and traces \cite{hu21trace}.

\begin{figure}[ht]
    \centering
    \scalebox{0.35}{\includegraphics{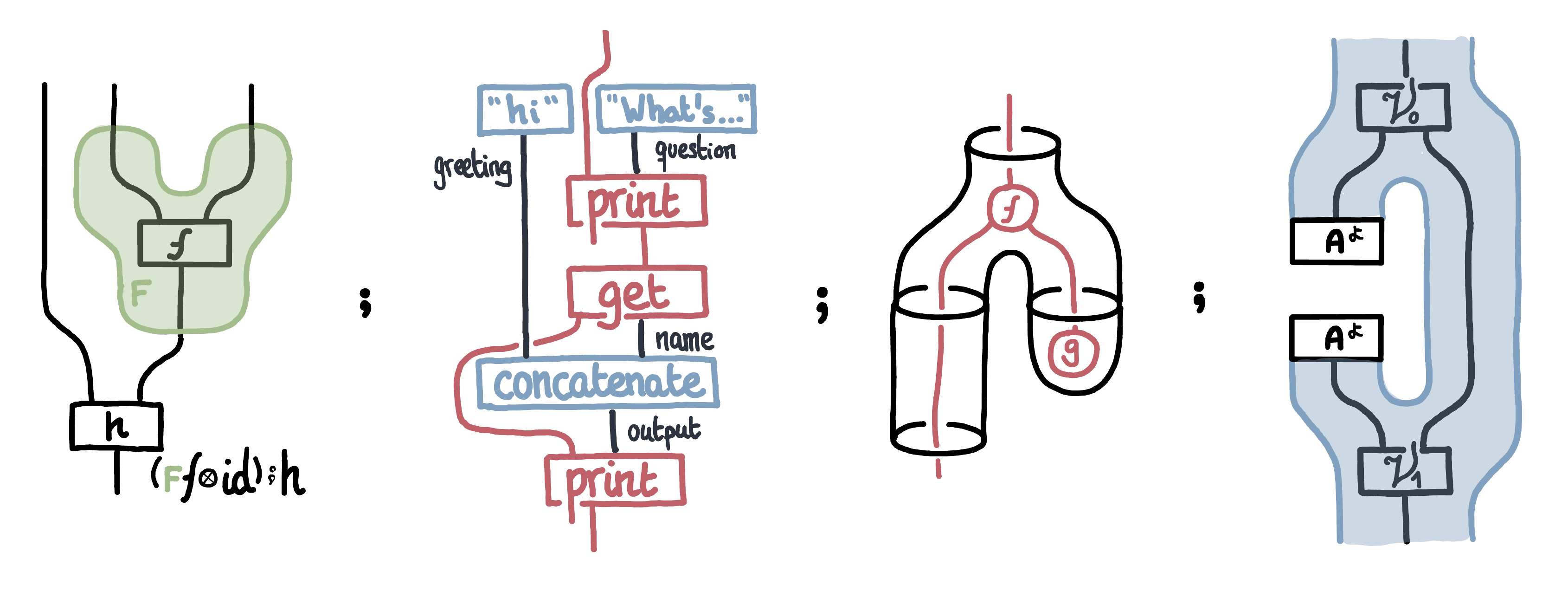}}
    \caption{Examples from the literature. From left to right: functor boxes \cite{mellies06}, premonoidal categories \cite{roman:effectful}, internal diagrams \cite{bartlett2015modular}, and combs or optics \cite{chiribella2009,clarke2022profunctor,hefford2022coend}.}
    \label{fig:examples}
\end{figure}

The extensions showcase the expressive power of string diagrams on surprisingly diverse  application domains.  At the same time, these different ideas could be regarded as separate ad-hoc extensions: they belong to different fields; they use different categorical formalisms. The overhead of learning and combining each one of them prevents the exchange of ideas between the different domains of application: e.g. an idea about topological quantum field diagrams does not transfer to premonoidal diagrams.

\paragraph{Collages.}
This manuscript claims that this division is only apparent and that all these extensions are particular instances of the same encompassing idea: that of glueing multiple string diagrams into what we call a \emph{collage of string diagrams}. We introduce a formal notion of collage (\Cref{sec:tricategory}) and employ string diagrammatic syntaxes for them, based on the calculus of bicategories (\Cref{sec:signaturebimodular,sec:functorboxsignature,sec:syntaxinternaldiags}).

Even though collages of string diagrams are our novel contribution, collages are not yet another new concept to category theory. ``Collage'' was Bob Walters' term for a lax colimit in a module-like category \cite{street1981cauchy}. This can be considered as a glueing of objects together along the action of a scalar. 
For example, given two sets $A$ and $B$, with an action of a monoid $M$, we can construct their tensor product $A \otimes_M B$, where $(a ⋅ m) ⊗ b = a ⊗ (m ⋅ b)$ for any scalar $m ∈ M$. Categorifying this idea in a possible direction we obtain monoidal categories acting on \emph{bimodular categories}. The following is the takeaway of this work.

\begingroup %
\addtolength\leftmargini{-0.05in}
\begin{leftbar} %
  \noindent Collages of string diagrams consist of multiple string diagrams of different monoidal categories glued together.
  Collages can be interpreted as \emph{pointed bimodular profunctors} between \emph{bimodular categories}.
\end{leftbar}
\endgroup

A \bimodularCategory{}, sometimes referred to as a biactegory \cite{capucci2022actegories}, is to a bimodule what a monoidal category is to a monoid. This is, a plain category $𝔸$ endowed with a left action of a monoidal category $(▹) \colon 𝕄 × 𝔸 → 𝔸$ and a right action of another, possibly different, monoidal category $(◃) \colon 𝔸 × ℕ → 𝔸$. We can \emph{collage} two bimodular categories along a common monoidal category that acts on both.
Later on the paper, exploiting a second axis of categorification, we pass from bimodular categories to \emph{bimodular profunctors}, which are a kind of 2-dimensional bimodule, and we define their collage. This structure facilitates glueing categories together in 2-dimensions: we can represent complexes of morphisms from different categories and glue them together. Collages of string diagrams are the syntactic representations of this glueing, in the same sense that ordinary string diagrams represent tensors in monoidal categories. 

We observe that collages of bimodular categories embed into a tricategory of pointed bimodules. This provides a versatile setting where we can interpret many syntaxes already present in the literature.

\paragraph*{Contributions.}
We introduce string diagrams of \bimodularCategories{} and we prove they construct the free bimodular category on a signature (\Cref{th:adj:bimod:bipoint}).
We introduce novel string diagrammatic syntax for \emph{functor boxes} and we prove it constructs the free lax monoidal functor on a suitable signature (\Cref{th:adjunctionfunctorboxes}).
We describe the tricategory of pointed bimodular profunctors (\Cref{def:tricategoryofpointedbimodularprofunctors}) and, in terms of it, we explain the semantics of functor boxes (\Cref{prop:bimodularlaxfunctor}) and internal diagrams (\Cref{th:semanticsopendiagrams}), for which we also provide a novel explicit formal syntax (\Cref{def:syntaxinternal}).

\section{String Diagrams of Bimodular Categories}
\label{sec:syntaxbimodularcats}

We introduce string diagrams for \bimodularCategories{} in terms of the better-known string diagrams of bicategories.
In algebra, a \emph{bimodule} is a structure with a compatible left and right action.
\emph{Bimodular categories are to bimodules what monoidal categories are to monoids}. 
Explicitly this means a category, $ℂ$, acted on by two monoidal categories, $𝕄$ and $ℕ$ \cite{skoda09:equivariant}.
Bimodular categories have also been known as ``biactegories'' \cite{capucci2022actegories,mccrudden00:categories}, while the name ``bimodule category'' typically refers to actions of certain vector enriched categories \cite{douglas19:balancedtensor}.
For our purposes, we consider a bimodular category $ℂ$, as gluing together the two acting categories, $𝕄$ and $ℕ$.

To simplify the presentation, we limit ourselves to considering only strict structure, but we expect that all of the results considered hold analogously in the weaker setting. In the following, we assume all monoidal, bimodular, and 2-categories to be strict, along with the associated functors between them.

\begin{definition}
  \defining{linkbimodularcategory}{}
  A \emph{bimodular category} $(ℂ,𝕄,ℕ)$ is a category $ℂ$ endowed with a left monoidal action 
  $(▹) \colon 𝕄 × ℂ \to ℂ$, and a right monoidal action $(◃) \colon ℂ × ℕ \to ℂ$, which are compatible in that
  $M ▹ (X ◃ N) = (M ▹ X) ◃ N$.

  \defining{linksbimod}{}
  \BimodularCategories{} over arbitrary monoidal categories form a category, $\sBimod$. The morphisms $(F, H, K) \colon (ℂ, 𝕄, ℕ) → (𝔻, ℙ, ℚ)$ consist of two monoidal functors $H \colon 𝕄 → ℙ$ and $K \colon ℕ → ℚ$ and a functor $F \colon ℂ → 𝔻$ that strictly preserves monoidal actions according to $H$ and $K$.
\end{definition}

Every monoidal category $(ℂ, \otimes, I)$ is a $(ℂ, ℂ)$-\bimodularCategory{} with its own tensor product defining the two actions.

\subsection{Signature of a Bimodular Category}
\label{sec:signaturebimodular}
The next sections exhibit a sound and complete string diagram syntax for bimodular categories. Bimodular string diagrams consist of two monoidal regions glued by a bimodular wire. We begin by defining a notion of bimodular signature and then construct an adjunction (\Cref{th:adjunctionconstruction}) using the notion of \emph{collages}.

\begin{definition}
  \label{def:bimodulargraph}\defining{linkbimodulargraph}{}
  A \emph{bimodular graph} $(𝓐,𝓜,𝓝)$ (the bimodular analogue of a multigraph \cite{shulman2016categorical}) is given by three sets of objects $(𝓐_{obj},𝓜_{obj},𝓝_{obj})$ and three different types of edges:
  \begin{itemize} 
    \item the left-acting edges, a set $𝓜(\vec{M} ; \vec{P})$ for each pair of lists of objects $\vec{M}, \vec{P} \in 𝓜_{obj}^*$,
    \item the right-acting edges, a set $𝓝(\vec{N} ; \vec{Q})$ for $\vec{N}, \vec{Q} \in 𝓝_{obj}^*$;
    \item the \emph{central edges}, a set of edges $𝓐(\vec{M}, A, \vec{N}\ ;\ \vec{P}, B, \vec{Q})$, for each $\vec{M}$, $\vec{P} \in 𝓜_{obj}^*$; each $\vec{N}$, $\vec{Q} \in 𝓝_{obj}^*$ and each $A$, $B ∈ 𝓐_{obj}$.
  \end{itemize}
\end{definition}

\begin{figure}[ht]
  \centering
  \scalebox{0.5}{\includegraphics{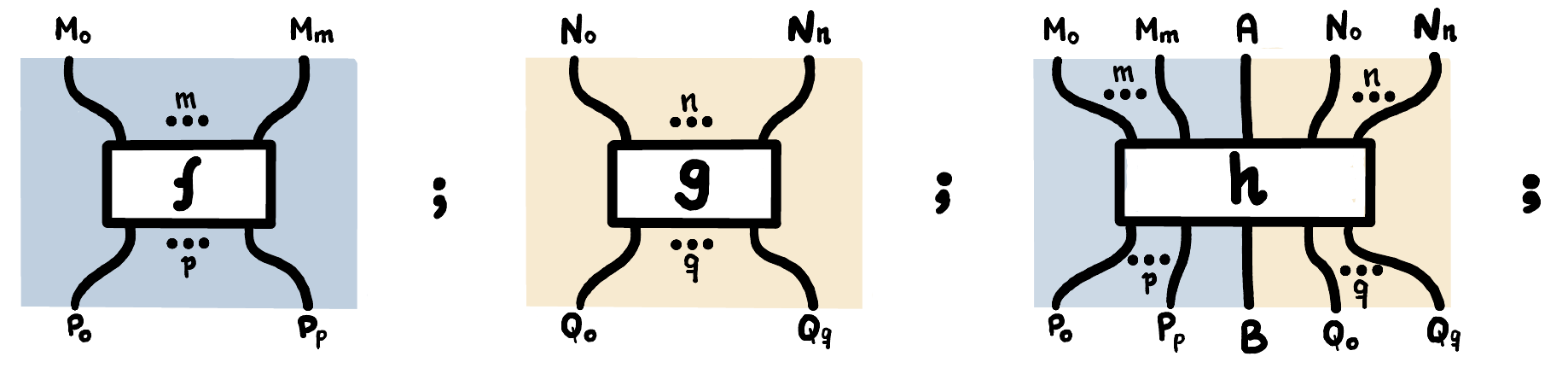}}
  \caption{Left, right, and central edges of a bimodular graph.}
  \label{fig:bimodulargraph}
\end{figure}

\begin{proposition}
  \label{prop:forgetfulbimodulargraph}
  \defining{linkbgraph}{}
  \BimodularGraphs{} form a category $\mathbf{bmGraph}$.
  We define a morphism of \bimodularGraphs{} 
    $(l,f,g) \colon (𝓐,𝓜,𝓝) \to (𝓐',𝓜',𝓝')$
  to be a triple of functions on objects,
    $(l_{obj},f_{obj},g_{obj})$,
  that extend to the morphism sets. 
  There exists a forgetful functor 
    $U \colon \sBimod → \bGraph$.
\end{proposition}

This provides a syntactic presentation of \bimodularCategories{}. We would like to additionally, construct a free model from a syntactic presentation.
We make use of the well-known similar result for 2-categories, exhibiting bimodular categories as certain bicategories: explicitly, those which are the \emph{collage} of a bimodular category.

\subsection{The Collage of a Bimodular Category}

Each profunctor induces a \emph{collage category}; analogously a \bimodularCategory{} induces a \emph{collage 2-category}.
This section proves that constructing the collage of a \bimodularCategory{} is left adjoint to considering the bimodular hom-category between any two cells of a 2-category.

\begin{definition}
  \label{def:collage}\defining{linkcollage}{}
  The \emph{collage} of an $(𝕄, ℕ)$-bimodular category $ℂ$ is a 2-category, $\Coll_ℂ$. This 2-category has two 0-cells, $M$ and $N$.
  The hom-categories are given by $\Coll_ℂ(M, M) = 𝕄$, $\Coll_ℂ(N, N) = ℕ$, and $\Coll_ℂ(M, N) = ℂ$; and finally $\Coll_ℂ(N, M)$ is the empty category.
  The composition of 1-cells is given by the monoidal products and actions.
\end{definition}

\begin{definition}
  \defining{linktwocattwo}{}
  The category of bipointed 2-categories, $\TwoCatTwo$, has as objects $(𝔸,M,N)$, 2-categories $𝔸$ with two chosen 0-cells on it, $M ∈ 𝔸$ and $N ∈ 𝔸$. 
  A morphism of $\TwoCatTwo$ is a 2-functor preserving the chosen 0-cells.
\end{definition}

\begin{theorem}
  \label{th:adj:bimod:bipoint}
  There exists an adjunction $\Coll_ℂ : \sBimod \rightleftarrows \TwoCatTwo : \Chosen$ given by the collage, and by picking the hom-category between the chosen 0-cells.
  Moreover, the unit of this adjunction is a natural isomorphism.
\end{theorem}

\subsection{String Diagrams of Bimodular Categories, via Collages}

We have the two ingredients for bimodular string diagrams:
sound complete string diagrams for 2-categories, and an embedding of bimodular categories into 2-categories by taking the collage.
We combine results to provide an adjunction from \bimodularGraphs{} to \bimodularCategories{}.

\begin{figure}[ht]
  \centering
\begin{tikzcd}
	\TwoGraph & {\TwoGraph_2} & \bGraph \\
	\TwoCat & \TwoCatTwo & \sBimod
	\arrow[""{name=0, anchor=center, inner sep=0}, "\Str"', curve={height=6pt}, from=1-1, to=2-1]
	\arrow[""{name=1, anchor=center, inner sep=0}, "{\mathsf{U}}"', curve={height=6pt}, from=2-1, to=1-1]
	\arrow[""{name=2, anchor=center, inner sep=0}, "{\Str_2}"', curve={height=6pt}, from=1-2, to=2-2]
	\arrow[""{name=3, anchor=center, inner sep=0}, "{\mathsf{U}_2}"', curve={height=6pt}, from=2-2, to=1-2]
	\arrow[""{name=4, anchor=center, inner sep=0}, "u"', curve={height=6pt}, tail, no head, from=1-3, to=1-2]
	\arrow[""{name=5, anchor=center, inner sep=0}, "i", curve={height=-6pt}, from=1-3, to=1-2]
	\arrow[""{name=6, anchor=center, inner sep=0}, "\Chosen", curve={height=-6pt}, from=2-2, to=2-3]
	\arrow[""{name=7, anchor=center, inner sep=0}, "\Coll", curve={height=-6pt}, from=2-3, to=2-2]
	\arrow[""{name=8, anchor=center, inner sep=0}, "\bmStr"', curve={height=6pt}, dashed, from=1-3, to=2-3]
	\arrow[""{name=9, anchor=center, inner sep=0}, "{\mathsf{U}}"', curve={height=6pt}, dashed, from=2-3, to=1-3]
	\arrow["\dashv"{anchor=center}, draw=none, from=0, to=1]
	\arrow["\dashv"{anchor=center}, draw=none, from=2, to=3]
	\arrow["\dashv"{anchor=center}, draw=none, from=8, to=9]
	\arrow["\dashv"{anchor=center, rotate=90}, draw=none, from=7, to=6]
	\arrow["\dashv"{anchor=center, rotate=90}, draw=none, from=5, to=4]
\end{tikzcd}
  \caption{{{Summary of adjunctions for the string diagrams of bimodular categories.}}}
  \label{fig:bimodularadjoints}
\end{figure}
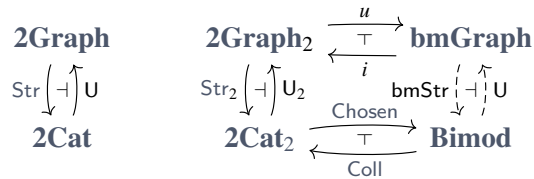
\begin{theorem}
  \label{th:adjunctionconstruction}
  There exists an adjunction between \bimodularGraphs{} and \bimodularCategories{}.
  The left adjoint finds the \bimodularCategory{} whose collage is the free 2-category on the \bimodularGraph{}, $\bmStr \colon \bGraph \to \sBimod$. The right adjoint is the forgetful functor  $\mathsf{U} \colon \sBimod → \bGraph$.
\end{theorem}

This result provides a basis for a graphical syntax for \bimodularCategories{}.
We now sketch an example of how these string diagrams can be of interest, but a larger class of examples come from premonoidal and effectful categories \cite{roman:effectful}.%

\subsection{Example: Shared State}

\begin{figure}[ht]
  \centering
  \scalebox{0.45}{\includegraphics{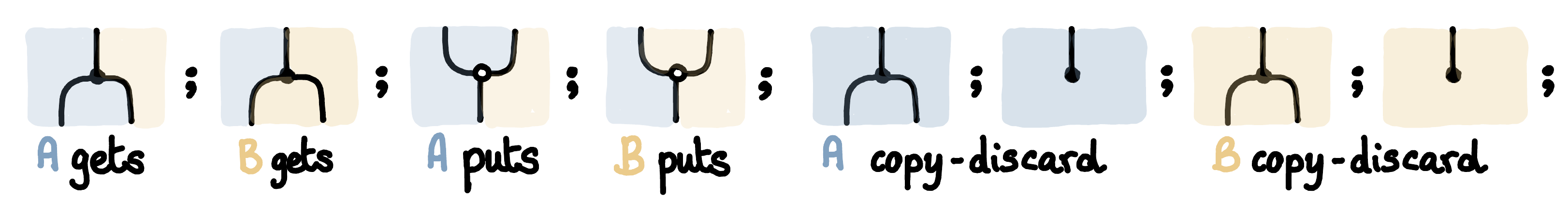}}
  \caption{{{Signature generators for the bimodular theory of shared state.}}}
  \label{fig:bimodularracesignature}
\end{figure}

In the same way that premonoidal categories are particularly well-suited to describe stateful computations, \bimodularCategories{} are particularly well-suited to describe shared state between two processes.
These processes can be different and even live in different categories.
As an example, consider the generators in \Cref{fig:bimodularracesignature}.
They represent two different process theories that access a common state with \textsf{get} and \textsf{put} operations. 

In the same way that monoidal categories are a good setting for defining monoids and comonoids, \bimodularCategories{} are a good setting for defining bimodules.
To capture interacting shared state, the generators of \Cref{fig:bimodularracesignature} are quotiented by the equations of a pair of semifrobenius modules with compatible comonoid actions and semimonoid actions.

\begin{figure}[ht]
  \centering
  \scalebox{0.35}{\includegraphics{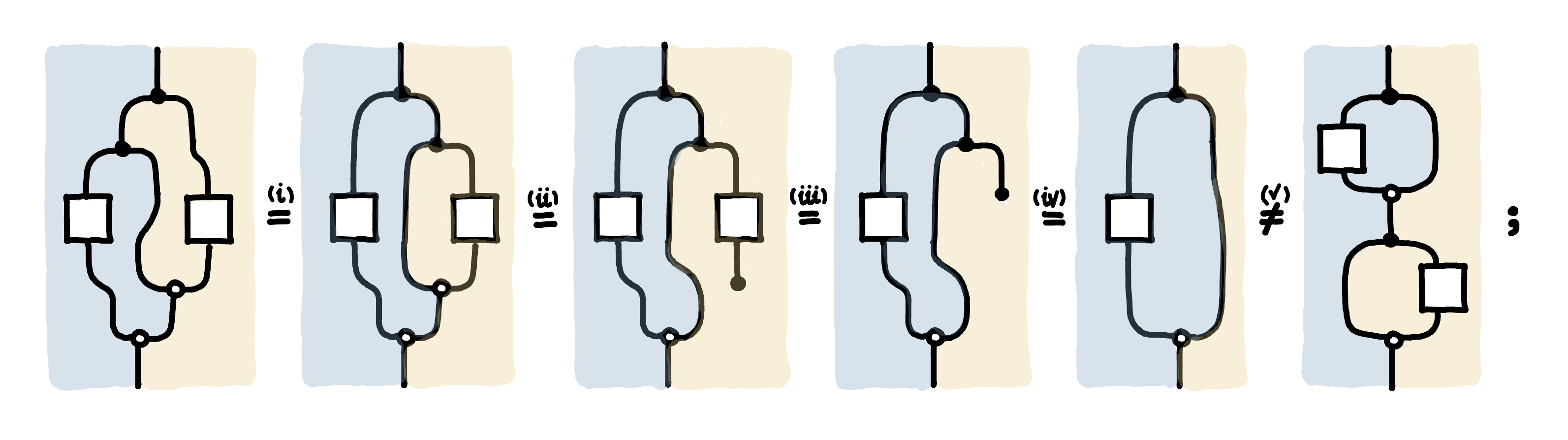}}
  \caption{{{Race condition in bimodular string diagrams.}}}
  \label{fig:bimodularrace}
\end{figure}
This setup is enough to exhibit one of the most salient features of shared state: \emph{race conditions}. Race conditions were first studied by Huffman in 1954, who used diagrams to show how the behavior of a shared state is dependent on the relative timing of the actions of the parties \cite{huffman1954synthesis}. We employ string diagrams of \bimodularCategories{} to show how two different timings of the actions -- the leftmost and rightmost sides of the equation in \Cref{fig:bimodularrace} -- result in two different executions: even when the two \textsf{get} statements are compatible \emph{(i)}, the two \textsf{put} statements interact causing the earlier of the two to be discarded \emph{(ii, iii, iv)}; this causes the discrepancy with the intended protocol \emph{(v)}.

\begin{figure}[ht]
  \centering
  \scalebox{0.45}{\includegraphics{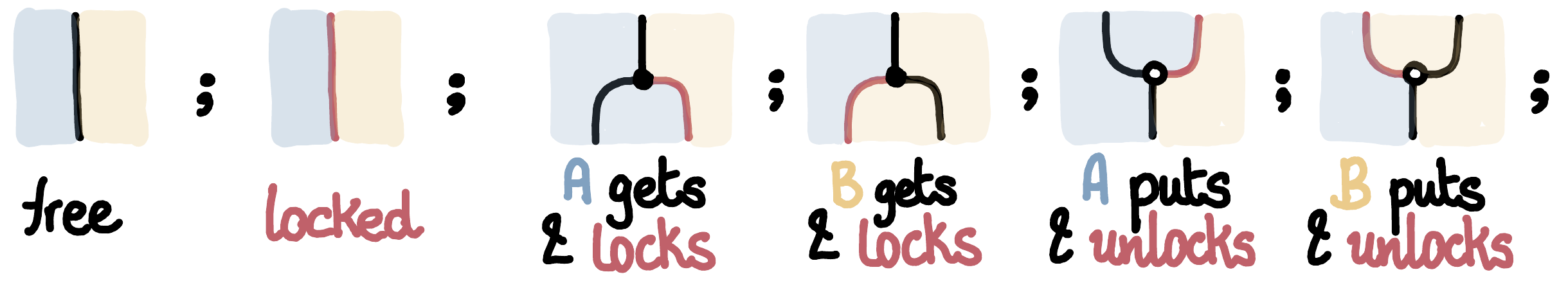}}
  \caption{{{Binary semaphore in bimodular string diagrams.}}}
  \label{fig:binarysemaphore}
\end{figure}
Race conditions have a commonly accepted workaround: the \emph{binary semaphore} \cite{SGG:operatingsystems}. Dijkstra described general semaphores with the aid of flow diagrams \cite{dijkstra1962over}; we instead use bimodular categories to model a binary semaphore (\Cref{fig:binarysemaphore}).
We consider a signature with two object generators, \textsf{free} and \textsf{locked}, for our bimodular category. 
Each operation must suitably lock or unlock the semaphore, rendering race conditions ill-typed, and leaving most of the interaction equations of the theory of shared state unnecessary.

String diagrams of bimodular categories model a pair of interacting monoidal categories.
We can also model an arbitrary number of interacting monoidal categories via a general collage construction.%
These collages can be subsumed into a `universe of collages' that we describe in \Cref{sec:bimodularprofunctors}: the tricategory of pointed bimodular profunctors. 
To motivate this, we study a second example: the syntax of functor boxes.
 
\section{String Diagrams of Functor Boxes}
\label{sec:syntaxfunctorbox}
Functor boxes are an extension of the string diagrammatic notation that represents plain functors, lax, oplax and strong monoidal functors. Functor boxes were introduced by Cockett and Seely \cite{cockett1999linearly} and later studied by Melliès \cite{mellies06}. We introduce here a syntactic presentation of (op)lax functor boxes that has the advantage of treating each piece of the box as a separate entity in a 2-category and applying the string diagrammatic calculus of 2-categories.

\subsection{Functor box signatures}
\label{sec:functorboxsignature}

\begin{definition}
  \defining{linkfbox}
  A \emph{functor box signature} $𝓕 = (𝓐,𝓧,𝓕_•,𝓕^•)$ consists of a pair of sets, $𝓐_{obj}$ and $𝓧_{obj}$, and four different types of edges:
  \begin{itemize}
    \setlength{\itemsep}{0em} 
    \item the plain edges, $𝓐(A_0,\dots,A_n ; B_0,\dots, B_m)$ for any objects $A_0,\dots,A_n,B_0,\dots,B_m ∈ 𝓐_{obj}$;
    \item the functor box edges, $𝓧(X₀, … ,Xₙ ; Y₀, … , Yₘ)$ for any objects $X_0,\dots,X_n,Y_0,\dots,Y_m ∈ 𝓧_{obj}$;
    \item the in-box edges, $𝓕_•(A₀, … ,Aₙ  ; Y₀, …, Yₘ)$ for any $A₀,…,Aₙ ∈ 𝓐_{obj}$ and $Y₀,…,Yₘ ∈ 𝓧_{obj}$ 
    \item the out-box edges, $𝓕^•(X₀, … ,Xₙ  ; B₀, …, Bₘ)$ for any $B₀,…,Bₘ ∈ 𝓐_{obj}$ and $X₀, …, Xₙ ∈ 𝓧_{obj}$.
  \end{itemize}
  A \emph{functor box signature morphism} $(h, k, l) \colon (𝓐,𝓧,𝓕) → (𝓑,𝓨,𝓖)$ is a pair of functions between the object sets, $h_{obj} \colon 𝓐_{obj} → 𝓑_{obj}$ and $k_{obj} \colon 𝓧_{obj} → 𝓨_{obj}$, that extend to a function between the edge sets;
  \begin{itemize}
    \item $h \colon 𝓐(A₀,…,Aₙ ; B₀,…,Bₘ) → 𝓑(h(A₀),…,h(Aₙ); h(B₀), …, h(Bₘ))$;
    \item $k \colon 𝓧(X₀,…,Xₙ ; Y₀,…,Yₘ) → 𝓨(k(X₀),…,k(Xₙ); k(Y₀), …, k(Yₘ))$;
    \item $l_• \colon 𝓕_•(A₀, … ,Aₙ  ; Y₀, …, Yₘ) → 𝓖_•(h(A₀), … ,h(Aₙ)  ; k(Y₀), …, k(Yₘ))$;
    \item $l^• \colon 𝓕^•(X₀, … ,Xₙ  ; B₀, …, Bₘ) → 𝓖^•(k(X₀), … ,k(Xₙ)  ; h(B₀), …, h(Bₘ))$.
  \end{itemize}
  Functor box signatures and homomorphisms form a category, $\mathbf{Fbox}$.
\end{definition}

\begin{figure}[ht]
  \centering
  \scalebox{0.45}{\includegraphics{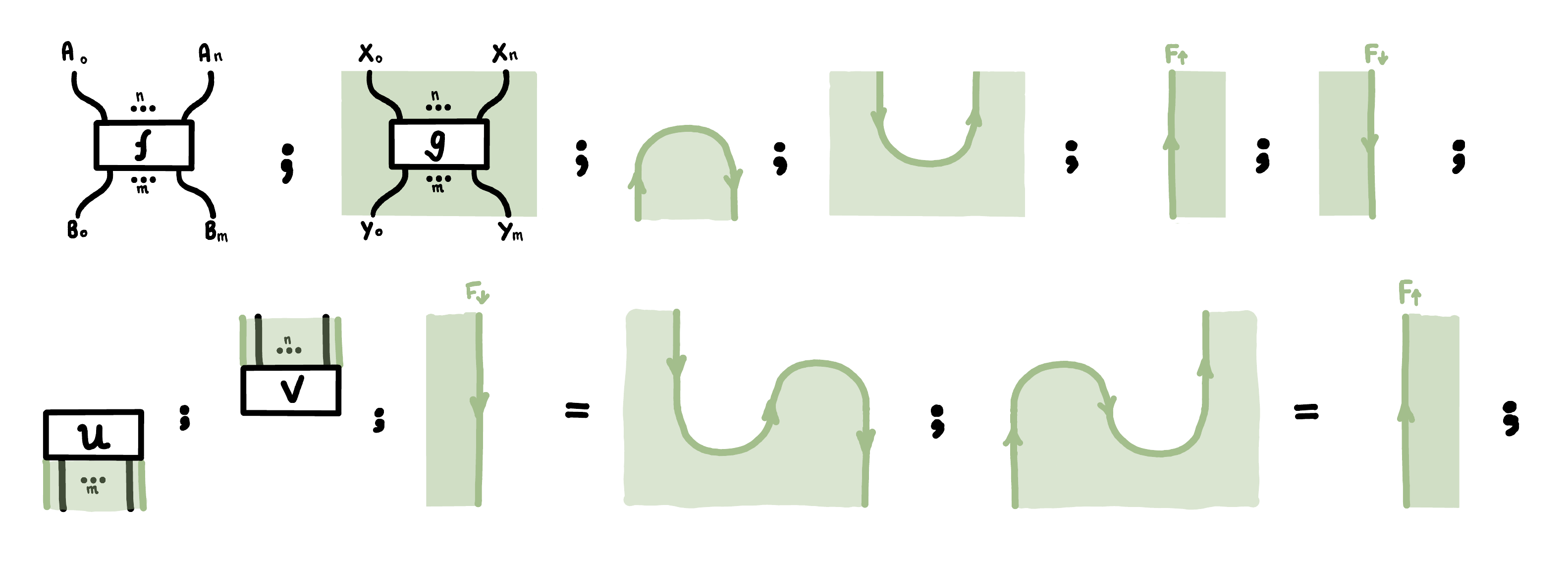}}
  \caption{Syntactic 2-category of a lax monoidal functor box signature.}
  \label{fig:functorboxsyntax}
\end{figure}

\begin{definition}
  The syntactic 2-category of a functor box signature $𝓕 = (𝓐,𝓧,𝓕_•,𝓕^•)$ is the 2-category freely presented by \Cref{fig:functorboxsyntax}, which we call $𝕊_{𝓐,𝓧,𝓕}$. 
  
  In other words, the 2-category $𝕊_{𝓐,𝓧,𝓕}$ contains exactly two 0-cells, labelled $𝓐$ and $𝓧$; it contains a 1-cell $A \colon 𝓐 → 𝓐$ for each $A ∈ 𝓐_{obj}$, a 1-cell $X \colon 𝓧 → 𝓧$ for each $X ∈ 𝓧_{obj}$ and, moreover, a pair of adjoint 1-cells $F^\uparrow \colon 𝓐 → 𝓧$ and $F^\downarrow \colon 𝓧 → 𝓐$. 
  Finally, it contains a pair of 2-cells witnessing the adjunction $F^\uparrow ⊣ F^\downarrow$, given by $n \colon \mathrm{id} → F^\uparrow ⨾ F^\downarrow$ and $e \colon F^\downarrow ⨾ F^\uparrow → \mathrm{id}$ which additionally satisfy the snake equations; and it also contains
  \begin{itemize}
    \item a 2-cell, $f ∈ 𝕊(𝓐,𝓐)(A₀ ⨾\dots ⨾ Aₙ ;\ B₀ ⨾ \dots ⨾ Bₘ)$, for each \emph{plain edge};
    \item a 2-cell, $g ∈ 𝕊(𝓧,𝓧)(X₀ ⨾ … ⨾ Xₙ ;\ Y₀ ⨾ … ⨾ Yₘ)$, for each \emph{functor box edge};
    \item a 2-cell, $u ∈ 𝕊(𝓐,𝓐)(A₀ ⨾ … ⨾ Aₙ ;\ F^\uparrow ⨾  Y₀ ⨾ … ⨾ Yₘ ⨾ F^\downarrow)$ for each \emph{in-box edge}; and
    \item a 2-cell, $v ∈ 𝕊(𝓐,𝓐)( F^\uparrow ⨾  X₀ ⨾ … ⨾ Xₙ ⨾ F^\downarrow ;\  B₀ ⨾ \dots ⨾ Bₘ)$ for each \emph{out-box edge}.
  \end{itemize}
\end{definition}

\subsection{Lax Monoidal Functor Semantics}

\begin{definition}[Lax functors category]
  \defining{linklaxcategory}
  An object of the \emph{lax functors category}, $\mathbf{Lax}$, is a pair of monoidal categories $(𝔸,𝕏)$ together with a lax monoidal functor between them, $(F,ε,μ)$; that is, a functor $F \colon 𝕏 → 𝔸$ endowed with two natural transformations
  $ε \colon I → FI, \mbox{ and } μ \colon FX ⊗ FY → F(X ⊗ Y),$ satisfying associativity $(μ ⊗ id) ⨾ μ = (id ⊗ μ) ⨾ μ$, left unitality $(ε ⊗ id) ⨾ μ = id$ and right unitality $(id ⊗ ε) ⨾ μ = id$.

  A morphism of the \emph{lax functors category}, from $(𝔸,𝕏,F,ε_F,μ_F)$ to $(𝔹,𝕐,G,ε_G,μ_G)$ is a pair of monoidal functors $H \colon 𝕏 → 𝔸$ and $K \colon 𝔸 → 𝔹$ such that $F ⨾ K = H ⨾ G$ and such that $K(ε_F) = ε_G$ and $K(μ_F) = μ_G$.
\end{definition}

\begin{theorem}
  \label{th:adjunctionfunctorboxes}
  There exists an adjunction between the category of functor box signatures, $\Fbox$, and the category of pairs of monoidal categories with a lax monoidal functor between them, $\Lax$. The free side of this adjunction is given by the syntax of \Cref{fig:functorboxsyntax}.
\end{theorem}

Collages, by themselves, explained the 2-region diagrams of \bimodularCategories{};
collages will also explain the two-region diagrams of functor boxes in \Cref{sec:functorboxesviacollages}.
However, as currently defined, collages are only sufficient to encode the vertical boundaries. 
To additionally represent boundaries along the horizontal axis we can make use of profunctors between bimodular categories and extend our notion of collage to these structures. %
Following this thread we find that collages embed into a tricategory of pointed bimodular profunctors, described in the next section, which we consider a universe of interpretation for all of the graphical theories described.

\section{Bimodular Profunctors}
\label{sec:bimodularprofunctors}
Where can we interpret all these string diagrams and provide compositional semantics for them? In this section, we introduce a single structure where all the previous calculi take semantics.

We will need two different ingredients: \emph{coends} and \emph{bimodularity}.
Coends and profunctors \cite{loregian2021,macLane1971}, far from being obscure concepts from category theory, can be seen as the right tool to glue together morphisms from different categories \cite{monoidalstreams,roman21}; we follow an explicitly \emph{pointed} version of coend calculus, which keeps track of the transformation between profunctors we are constructing (\Cref{sec:pointedprofunctors}).
In a similar sense, \bimodularCategories{} tensor together objects from different monoidal categories.
Both ideas combine into the calculus of pointed bimodular profunctors.

\subsection{Bimodular Profunctors}

Consider $ℂ$ and $𝔻$, both $(𝕄,ℕ)$-bimodular categories.
A natural notion of morphism between them is a functor $ℂ \to 𝔻 $ which preserves both actions.
However, there is another notion of morphism between them, which is a generalization of a profunctor between categories to this bimodular setting.
\BimodularProfunctors{} are a generalized reformulation of the Tambara modules of Pastro and Street \cite{pastro07doubles}.

\begin{definition}
  \label{def:bimodularprofunctor}\defining{linkbimodularprofunctor}{}
  Let $𝕄$ and $ℕ$ be two monoidal categories and
  let $ℂ$ and $𝔻$ be two $(𝕄,ℕ)$-bimodular categories.
  A \emph{bimodular profunctor} from $ℂ$ to $𝔻$ is a profunctor $T \colon ℂ\op \times 𝔻 \to \Set$ with a natural family of strengths,
    $$t_M : T(X, Y) \to T(M ▹ X, M ▹ Y),\quad\mbox{ and }\quad t^N : T(X, Y) \to T(X ◃ N, Y ◃ N),$$
  such that the actions are associative, $t_M ⨾ t_{M'} = t_{M ⊗ M'}$ and $t_N ⨾ t_{N'} = t_{N ⊗ N'}$, unital $t_I = id$ and $t^I = id$, and compatible, $t_M ⨾ t^N = t^N ⨾ t_M$, up to the coherence isomorphisms of the monoidal category. 
\end{definition}

\begin{proposition}
  For any pair of monoidal categories, $𝕄$ and $ℕ$, there is a 2-category $\Mod{𝕄}{ℕ}$ of $(𝕄, ℕ)$-\bimodularCategories{}, \bimodularProfunctors{}, and natural transformations between them. 
\end{proposition}
  
  These will form the hom-bicategories of the tricategory we later define.
  The other significant piece of data we require is a family of tensors $\otimes : \Mod{𝕄}{ℕ} \times \Mod{ℕ}{𝕆} \to \Mod{𝕄}{𝕆}$, which we now study. 
\subsection{Tensor of Bimodular Profunctors}

The tensor of \bimodularCategories{} is similar to the tensor of modules over a monoid in classical algebra:
we consider pairs of elements and we quotient out the action of a common scalar \cite{ostrik03:modulecategories}. 
In this case, the quotienting is substituted by an appropriate structural isomorphism: the \emph{equilibrator}.

\begin{definition}[Tensor of bimodular categories, {{\cite{ostrik03:modulecategories}}}]
  Let $ℂ$ be a $(𝕄,ℕ)$-\bimodularCategory{} and let $𝔻$ be a $(ℕ,𝕆)$-\bimodularCategory{}.
  Their tensor product, $ℂ \otimes_ℕ 𝔻$, is a category with the same objects as $ℂ × 𝔻$: we write them as $X \otimes_ℕ Y$. The category is presented by the morphisms of $ℂ × 𝔻$ and a free family of natural isomorphisms, called the \emph{equilibrators},
  $$τ_{X,N,Y} \colon (X \triangleleft N) \otimes_N Y \to X \otimes_N (N \triangleright Y), \mbox{ for each } N \in ℕ, X \in ℂ, Y \in 𝔻,$$
  which are additionally quotiented by the following equations up to the structure isomorphisms of the monoidal actions, $τ_{X, M \otimes N, Y} = τ_{X \triangleleft M, N, Y} ⨾ τ_{X, M, N \triangleright Y}$, and $τ_{X, I, Y} = \mathrm{id}$.
\end{definition}

\begin{definition}
  Let $ℂ$ and $ℂ'$ be two $(𝕄,ℕ)$-\bimodularCategories{} and let $𝔻$ and $𝔻'$ be a $(ℕ,𝕆)$-\bimodularCategories{}.
  Given two \bimodularProfunctors{}, $T \colon ℂ \to ℂ'$ and $R \colon 𝔻 \to 𝔻'$, their tensor is a \bimodularProfunctor{}, $T \otimes_ℕ R \colon ℂ \otimes_ℕ 𝔻 \to ℂ' \otimes_ℕ 𝔻'$, defined by
  $$(T \otimes_ℕ R)(X \otimes_N Y;X' \otimes_N Y') = T(X;X') × R(Y,Y') / (\sim),$$ 
  where $(\sim)$ is the equivalence relation generated by $(t_N(x),y) \sim (x,t_N(y))$.
\end{definition}

\subsection{Pointed Profunctors}
\label{sec:pointedprofunctors}

Profunctors deal with families of morphisms, and their natural isomorphisms determine correspondences between these families.
However, when we use profunctors for the semantics of string diagrams, we most often want to single out a particular morphism between a particular pair of objects.
A simple technique to achieve this is to use \emph{pointed profunctors} instead of simply profunctors: this technique was explicitly described by this second author \cite{roman21} although it has implicit appearances in the literature \cite{bartlett2015modular,hu21trace}.

\begin{definition}
  A pointed profunctor $(P,p) \colon (𝔸, X) \to (𝔹, Y)$ between two pointed categories with a chosen object $X \in 𝔸_{obj}$ and $Y \in 𝔹_{obj}$ is a profunctor $P \colon 𝔸 \to 𝔹$ together with an element $p \in P(A, B)$ of the profunctor evaluated on the chosen object of the categories.
\end{definition}

From now on, we work using pointed profunctors instead of plain profunctors. %

\subsection{The Tricategory of Pointed Bimodular Profunctors}
\label{sec:tricategory}
We call \emph{collages of string diagrams} to the diagrams of the tricategory of pointed bimodular profunctors.

\begin{definition}
  \label{def:tricategoryofpointedbimodularprofunctors}
  The tricategory of pointed bimodular profunctors, $\mathbb{BmProf}\pointed$, has as 0-cells the monoidal categories, $𝕄, ℕ, 𝕆,\dots$.
  The 1-cells between two monoidal categories $𝕄$ and $ℕ$ are \emph{pointed bimodular categories}, $(𝔸,▹,◃,A)$, consisting of a $(𝕄,ℕ)$-bimodular category with two actions $(𝔸,▹,◃)$ and some object of that category, $A ∈ 𝔸$.
  Pointed bimodular categories compose by the tensor of bimodular categories,
  $$(𝔸,▹,◃,A) ⊗_{ℕ} (𝔹,▹,◃,B) = (𝔸 ⊗_{ℕ} 𝔹, ▹,◃, A ⊗_{ℕ} B).$$
  The 2-cells between two pointed bimodular categories $(𝔸,▹,◃,A)$ and $(𝔹,▹,◃,B)$ are \emph{pointed bimodular profunctors} $(P,t,p)$, consisting of a profunctor $P \colon 𝔸 → 𝔹$ together with a point $p ∈ P(A,B)$ that are moreover bimodular with compatible natural transformations
    $t_M \colon P(A;B) \to P(M ▹ A; M ▹ B)$, and 
    $t_N \colon P(A;B) \to P(A ◃ N; B ◃ N)$.
  These 2-cells compose by profunctor composition and by the tensor of bimodular profunctors.
  
  Finally, the 3-cells between two pointed bimodular profunctors $(P,t,p)$ and $(Q,r,q)$ are bimodular natural transformations that preserve the point, consisting of a natural transformation $α \colon P → Q$ such that the $α(p) = q$ and, moreover, $t_M ⨾ α = α ⨾ r_M$ and $t_N ⨾ α = α ⨾ r_N$.
\end{definition}

\begin{remark}
At the moment of writing, it is unclear to the authors whether a string diagrammatic calculus for tricategories, described by transformations of the string diagrammatic calculus of bicategories, has been fully described and proved sound and complete. However, there seems to be a consensus that this would be the right language for tricategories: much literature assumes it. Let us close this section by tracking explicitly the assumptions we need to employ a diagrammatic syntax for bimodular profunctors.
\end{remark}

\begin{conjecture}
  The previous data satisfies all coherence conditions of a tricategory.
  Moreover, we can reason with tricategories using the calculus of deformations of string diagrams, extending the string diagrams for quasistrict monoidal 2-categories of Bartlett \cite{bartlett14}.
\end{conjecture}

\subsection{Functor Boxes via Collages of String Diagrams}
\label{sec:functorboxesviacollages}
\begin{figure}[ht]
  \centering
  \scalebox{0.30}{\includegraphics{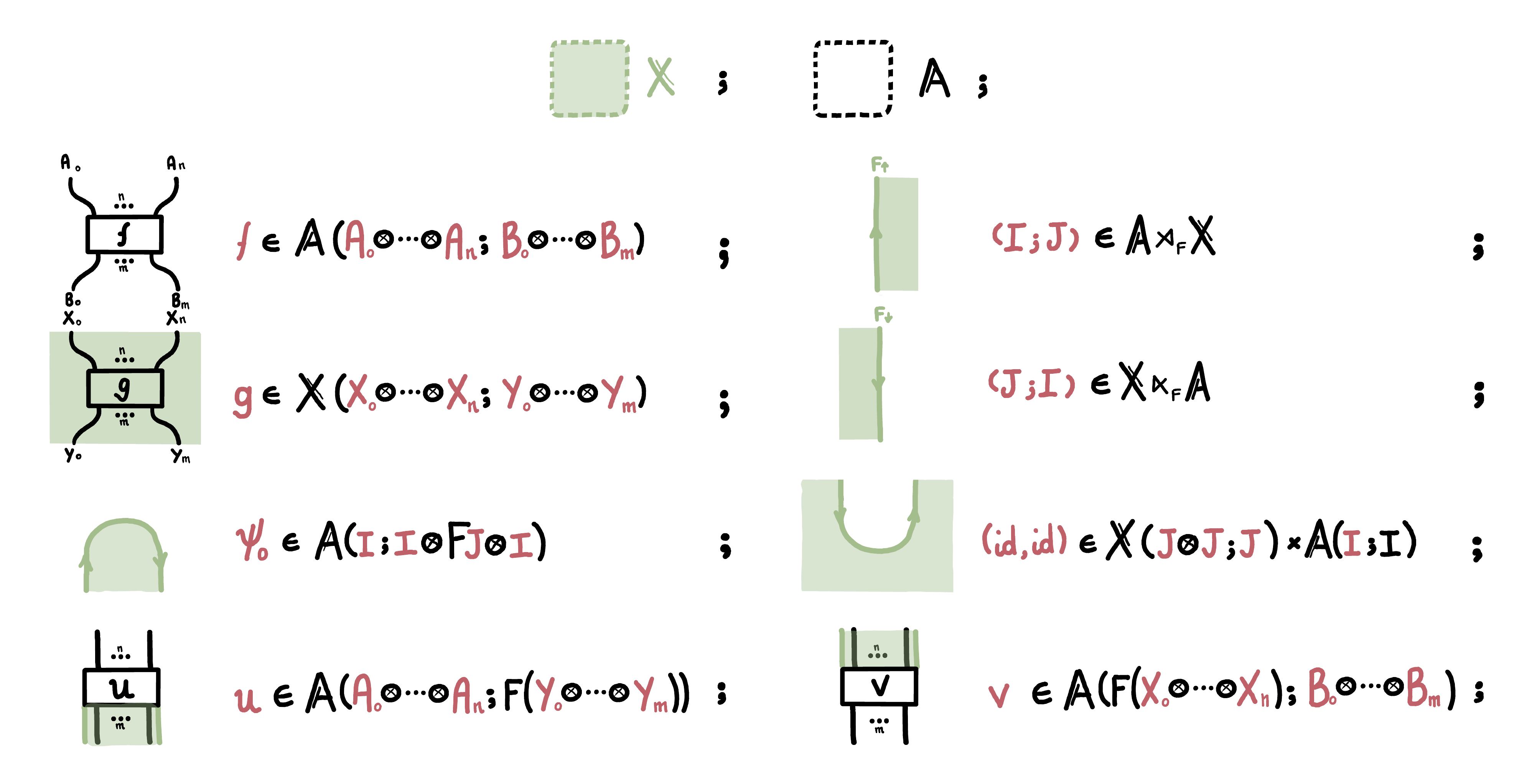}}
  \caption{Semantics for functor boxes in terms of pointed bimodular profunctors.}
  \label{fig:functorboxsemantics}
\end{figure}

The following \Cref{fig:functorboxsemantics} details how to interpret functor boxes as collages of string diagrams.
The colored region represents the domain of the lax monoidal functor; the white region represents the codomain.
Morphisms of both categories are interpreted as elements of their respective hom-profunctors, and the laxators are used to merge colored regions. The only element that we will explicitly detail is the \bimodularCategory{} that appears in the closing and opening wires of a functor box.

\begin{proposition}[Bimodular categories of a lax monoidal functor]
  \label{prop:bimodularlaxfunctor}
  Let $\mathbb{X}$ and $𝔸$ be two monoidal categories and let $F \colon \mathbb{X} \to 𝔸$ be a monoidal functor between them, endowed with natural transformations $\psi_0 \colon J \to FI$ and $\psi_2 \colon FX \otimes FY \to F(X \otimes Y)$.
  The following profunctors, $𝔸 \mathbin{\rtimes_F} \mathbb{X} \colon 𝔸 \times \mathbb{X} \to 𝔸 \times \mathbb{X}$ and $\mathbb{X} \mathbin{\ltimes_F} 𝔸 \colon \mathbb{X} \times 𝔸 \to \mathbb{X} \times 𝔸$ determine two promonads, and therefore two Kleisli categories.
  \begin{align*}
    𝔸 \mathbin{\rtimes_F} \mathbb{X} (A,X;B,Y) = \int^{M \in \mathbb{X}} 𝔸(A ; B \otimes FM) \times \mathbb{X}(M \otimes X ; Y); \\
    \mathbb{X} \mathbin{\ltimes_F} 𝔸 (X,A;Y,B) = \int^{M \in \mathbb{X}} 𝔸(A ; FM \otimes B) \times \mathbb{X}(M \otimes A ; B); 
  \end{align*}
  These two Kleisli categories are $(𝔸,\mathbb{X})$ and $(\mathbb{X},𝔸)$-bimodular, respectively.
\end{proposition}

\section{String Diagrams of Internal Diagrams}
\label{sec:syntaxinternaldiags}
The tubular 3-dimensional cobordisms of internal diagrams are first described as a Frobenius algebra by Bartlett, Douglas, Schommer-Pries and Vicary \cite{bartlett2015modular}. We are indebted to this first introduction, which made internal diagrams into a convenient graphical notation in topological quantum field theory \cite{bartlett2015modular}. Internal diagrams themselves were later given explicit semantics in a monoidal bicategory of pointed profunctors; this was the subject of this second author's contribution to \emph{Applied Category Theory 2020} \cite{roman2020}. 
An important aspect of the syntax of internal diagrams is their 3-dimensional nature: the syntax not only contains string diagrams but also reductions between them. 

\begin{figure}[ht]
  \centering
  \scalebox{0.4}{\includegraphics{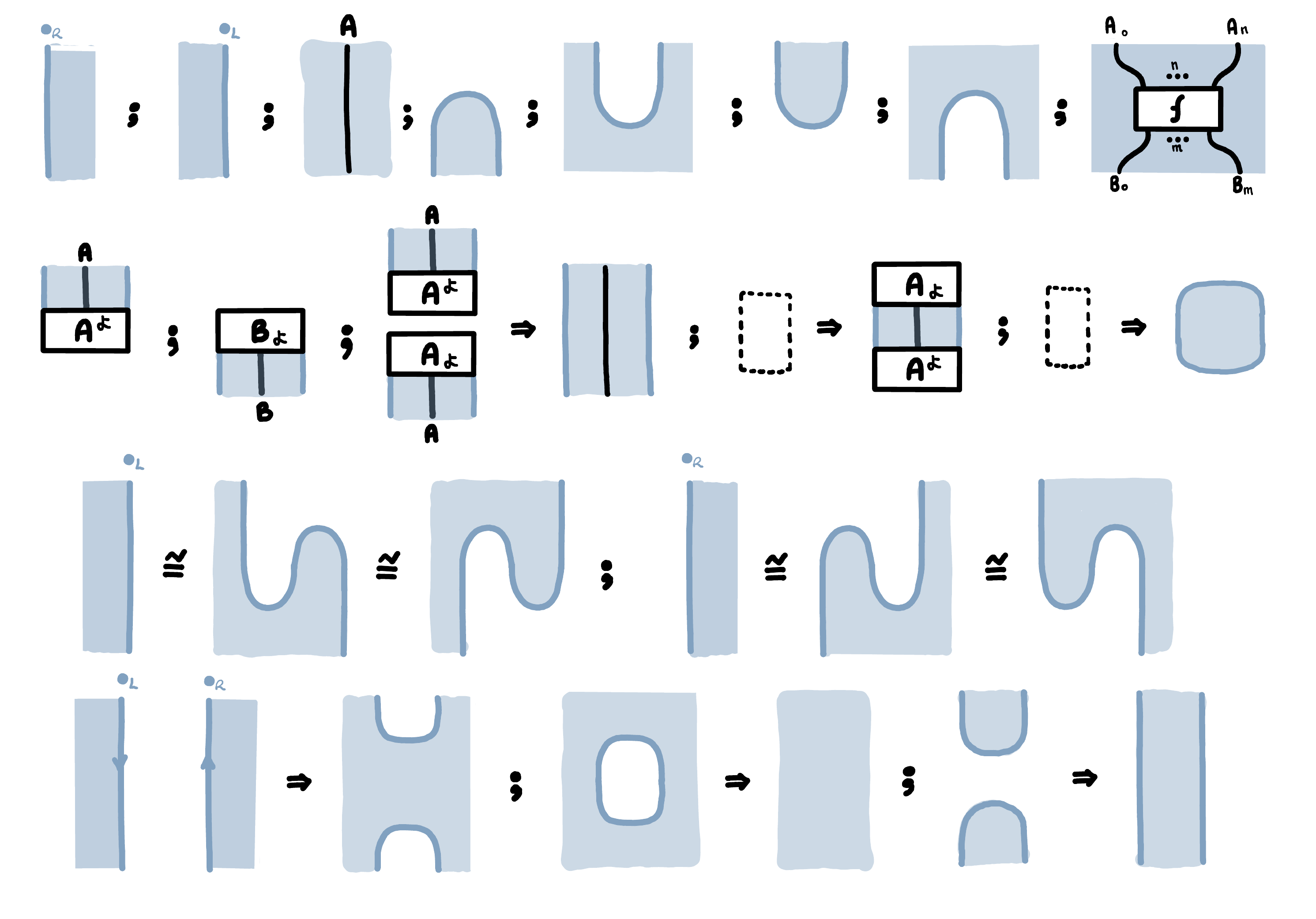}}
  \caption{Syntax for open internal diagrams.}
  \label{fig:opendiagramssyntax}
\end{figure}

We introduce here a novel syntactic presentation of \emph{internal diagrams} that has the advantage of treating each piece of an internal diagram (including the closing and opening of tubes) as a separate entity in a tricategory. That is, the identity tube or the multiplication and comultiplication tubes are constructed out of smaller pieces in \Cref{fig:opendiagramssyntax}. As a consequence, we are later able to introduce, for the first time, a more refined semantics in terms of a tricategory of \emph{pointed bimodular profunctors}.

\begin{definition}
  \defining{linkpolygraph}{}
  A \emph{polygraph}, $𝓖$, is the signature for the string diagrams of a monoidal category. %
  It consists of a set of objects, $𝓖_{obj}$, and a set of morphisms $𝓖(A₀, …, Aₙ; B₀, …, Bₘ)$ between any two lists of objects, $A₀,…,Aₙ,B₀,…,Bₘ ∈ 𝓖_{obj}$.
\end{definition}

\begin{definition}
  \label{def:syntaxinternal}
  The \emph{syntactic 3-category of internal diagrams} over a \polygraph{} $𝓖$ is the 3-category $\mathbf{G}$ presented by the cells in \Cref{fig:opendiagramssyntax}. In other words, it contains two 0-cells, $𝓘$ and $𝓖$, in white and blue in the figure, respectively. It contains a 1-cell $A \colon 𝓖 → 𝓖$ for each object $A ∈ 𝓖_{obj}$ and two 1-cells, $L_• \colon 𝓘 → 𝓖$ and $R_• \colon 𝓖 → 𝓘$ forming two 2-adjunctions $(L_•) ⊣ (R_•)$ and $(R_•) ⊣ (L_•)$ up to a 3-cell. It contains the following 2-cells,
  \begin{itemize} 
    \item two 2-cells $n₁ \colon \mathrm{id} → L_• ⨾ R_•$ and $e₁ \colon R_• ⨾ L_• → \mathrm{id}$ witnessing the 2-adjunction $(L_•) ⊣ (R_•)$ and two 2-cells $n₂ \colon 1 → R_• ⨾ L_•$ and $e₂ \colon  L_• ⨾ R_• → \mathrm{id}$ witnessing the 2-adjunction $(R_•) ⊣ (L_•)$ -- see Vicary and Heunen \cite{heunen2019categories} for a reference on 2-adjunctions and the swallowtail equations;
    \item two 2-cells, $A^{ょ} \colon L_• ⨾ A ⨾ R_• → \mathrm{id}$ and $A_{ょ} \colon \mathrm{id} → L_• ⨾ A ⨾ R_•$, forming an adjunction $A^{ょ} ⊣ A_{ょ}$ for each object $A ∈ 𝓖_{obj}$; and a 2-cell, $f \colon A₀ ⨾ … ⨾ Aₙ → B₀ ⨾ … ⨾ Bₘ$, for each edge $f ∈ 𝓖(A₀,…,Aₙ;B₀,…,Bₘ)$.
  \end{itemize}
  Finally, it contains the following 3-cells,
  \begin{itemize}
    \item two invertible 3-cells, $\alpha₁ \colon (1 \otimes n₁) ⨾ (e₁ \otimes 1) \to 1$ and $\beta₁ \colon (n₁ \otimes 1) ⨾ (1 \otimes e₁) \to 1$, witnessing the 2-adjunction $(L_•) ⊣ (R_•)$ and satisfying the swallowtail equations; and two invertible 3-cells, $\alpha'₂ \colon (1 \otimes n_2) ⨾ (e₂ \otimes 1) \to 1$ and $\beta₂ \colon (n₁ \otimes 1) ⨾ (1 \otimes e₁) \to 1$, witnessing the 2-adjunction $(R_•) ⊣ (L_•)$ and and satisfying the swallowtail equations;
    \item two 3-cells, $c \colon A^{ょ} ⨾ A_{ょ} \to 1$ and $i \colon 1 \to A_{ょ} ⨾ A^{ょ}$, witnessing the adjunction $A^{ょ} ⊣ A_{ょ}$ and satisfying the snake equations;
    \item two 3-cells, $uᵢ \colon n₁ ⨾ e₂ \to 1$ and $v_i \colon 1 \to e₂ ⨾ n₁$ witnessing an adjunction $e_2 \dashv n₁$ and satisfying the snake equations; two 3-cells $uⱼ \colon 1 \to n₂ ⨾ e₁$ and $v_i \colon e₁ ⨾ n₂ \to 1$ witnessing an adjunction $n₂ \dashv e₁$ and satisfying the snake equations.  
  \end{itemize}
\end{definition}

\begin{theorem}
  \label{th:semanticsopendiagrams}
  For any interpretation of a polygraph into a monoidal category, there exists a 3-functor from the syntactic tricategory of internal diagrams into pointed bimodular profunctors that preserves this interpretation.
\end{theorem}

\begin{remark}
This syntax can be exemplified by evaluating a quantum comb \cite{chiribella2009}, or a monoidal lens \cite{riley2018categories} with a morphism, in terms of internal string diagrams \cite{hu21trace}, see \Cref{fig:combs}.  It has been used more generally to reason about coends in monoidal categories \cite{roman21} and topological quantum field theory \cite{bartlett2015modular}.
  
\begin{figure}[ht]
  \centering
  \scalebox{0.3}{\includegraphics{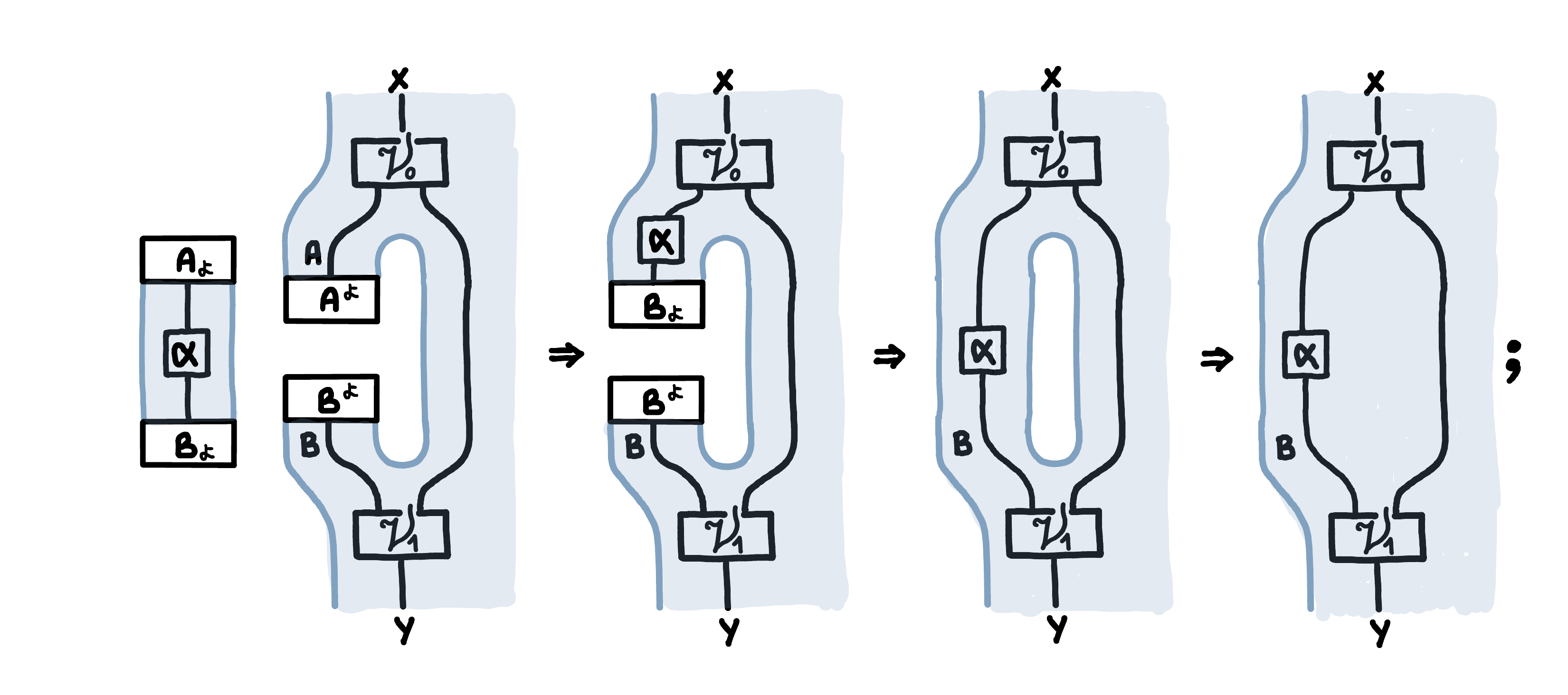}}
  \caption{Evaluating a comb in terms of internal string diagrams.}
  \label{fig:combs}
\end{figure}
\end{remark}

\section{Conclusions}

Collages of string diagrams provide an abundant graphical calculus. Functor boxes, tensors of bimodular categories and internal diagrams all exist in the graphical calculus of collages.
Their technical underpinning is complex: we characterized them as diagrams of pointed bimodular profunctors, but these arrange themselves into a tricategory, which may be difficult to reason about.

Apart from introducing the technique of collages and formalizing multiple extensions to string diagrams, we would like to call attention to the techniques we use: most of our results on soundness and completeness of diagrams are arranged into adjunctions, which allows us to prove them by reusing the better-known results on soundness and completeness for monoidal categories and bicategories.

\paragraph{Related work.}
An important line of research revolves around \emph{module categories} and \emph{fusion categories}, some specific enriched categories with actions with applications in topological quantum field theories \cite{douglas19:balancedtensor,drinfeld10,ostrik03:modulecategories}. Especially relevant and recent is Hoek's work, which constructs diagrams for a bimodule category \cite[Theorem 3.5.2]{hoek2019drinfeld}.
We follow the more elementary notion of \bimodularCategory{}, called ``\emph{biactegory}'' in the taxonomy of Capucci and Gavranovi\'c \cite{capucci2022actegories}. 
Cockett and Pastro \cite{cockettP09} have used instead \emph{linear actions} for concurrency, and even when we take inspiration from their work, their approach is more sophisticated and expressive than our toy example demonstrating bimodular categories (\Cref{fig:bimodularrace}).

Most work has been presented for some particular cases of collages: functor boxes have been extensively employed, but never reduced to string diagrams \cite{cockett1999linearly,mellies06}; internal diagrams have served both quantum theory and category theory \cite{bartlett2015modular,hu21trace,lobskizanasi:layered}, and can be given semantics into pointed profunctors \cite{roman2020}, but again a presentation as string diagrams was missing. A convenient algebra of lenses \cite{riley2018categories}, a particular type of incomplete diagram, has been recently introduced \cite{earnshaw2023produoidal}, but this is still independent of the semantics of arbitrary internal diagrams.

Finally, the first author has published a blog post that accompanies this manuscript \cite{braithwaite2023}.

\paragraph{Further work.}
It should be possible to ``destrictify'' many of the results of this paper. 
We have only presented a 1-adjunction between strict \bimodularCategories{} and bipointed 2-categories, but a higher adjunction would allow us to reuse coherence for bicategories to automatically obtain coherence for \bimodularCategories{}.
We indicated along the paper the conjectures where further work is warranted.

We conjecture that pointed bimodular profunctors form a compact closed tricategory, with the dual of each monoidal category being the \emph{reverse monoidal category}, $A \otimes_{Rev} B = B \otimes A$.
Even when it may be conceptually clear what a compact tricategory should be, it is technically challenging to come up with a concrete definition for it in terms of coherence equations.

\section*{Acknowledgements}
The authors want to thank David A. Dalrymple for discussion on the string diagrammatic interpretation of functor boxes; and Matteo Capucci for several insightful conversations about notions of 2-dimensional profunctor, that helped us understand how to tie disparate aspects of this story together.
The authors thank John Baez, the editors, and the anonymous reviewers at ACT23 for multiple comments and suggestions that improved this manuscript.

Dylan Braithwaite was supported by an Industrial CASE studentship from the UK Engineering and Physical Sciences Research Council (EPSRC) and the National Physical Laboratory. 
Mario Román was supported by the European Union through the ESF Estonian IT Academy research measure (2014-2020.4.05.19-0001).

\newpage
\bibliographystyle{eptcs}
\bibliography{bibliography}

\end{document}